\def\ZZ         {{\bf Z}}

\def\CC         {{\bf C}}
\def\QQ         {{\bf Q}}
\def\PP         {{\bf P}}
\def\ii         {{\rm i}}
\def\ee         {{\rm e}}
\def\tr         {{\rm tr}}

\def\dim        {{\rm dim}}
\def\sin        {{\rm sin}}

\def\deg        {{\rm deg}}
\def\codim        {{\rm codim}}
\def\Ell        {{\cal ELL}}

\def\cal        {\mathcal}

\documentclass{amsart}
\usepackage{amssymb}

\newtheorem{theorem}{Theorem}[section]
\newtheorem{lemma}[theorem]{Lemma}
\newtheorem{proposition}[theorem]{Proposition}
\newtheorem{corollary}[theorem]{Corollary}

\theoremstyle{definition}
\newtheorem{definition}[theorem]{Definition}
\newtheorem{conjecture}[theorem]{Conjecture}

\theoremstyle{remark}
\newtheorem{remark}[theorem]{Remark}

\numberwithin{equation}{section}

\begin{document}

\title{Elliptic Genera of Singular Varieties}

\author{Lev Borisov}
\address{Department of Mathematics\\
Columbia University\\
New York, NY  10027}
\email{lborisov@math.columbia.edu}

\author{Anatoly Libgober}
\address{Department of Mathematics\\
University of Illinois\\
Chicago, IL 60607}
\email{libgober@math.uic.edu}

\subjclass{to be included}

\date{to be included}


\begin{abstract}
{Orbifold elliptic genus and elliptic genus of singular varieties 
 are introduced and relation between them is studied. 
 Elliptic genus of singular varieties is given in terms 
 of a resolution of singularities and extends the elliptic 
genus of Calabi-Yau hypersurfaces in Fano Gorenstein toric
varieties introduced earlier. Orbifold elliptic genus is 
given in terms of the fixed point sets of the action.
We show that the generating function for this orbifold 
elliptic genus  $\sum Ell_{orb}(X^n,\Sigma_n)p^n$ for symmetric groups 
$\Sigma_n$ acting on $n$-fold products coincides with 
the one proposed by Dijkgraaf, Moore, Verlinde and Verlinde. 
Two notions of elliptic genera are conjectured to coincide.
}
\end{abstract}

\maketitle

\section{Introduction}
This work started as an attempt to 
understand a beautiful formula for the generating 
function for the orbifold elliptic genera of symmetric 
products due to R.Dijkgraaf, G.Moore, E.Verlinde and H.Verlinde,
which looks as follows (cf. \cite{DMVV}):
\begin{equation}
\sum_{n\geq 0}p^n Ell_{orb}(X^n/\Sigma_n;y,q)=
\prod_{i=1}^\infty \prod_{l,m}
\frac 1
{(1-p^iy^lq^m)^{c(mi,l)}}.
\label{dmvv}
\end{equation}
Here $X$ is a K\"{a}hler manifold, $\Sigma_n$ is the symmetric group acting on the
$n$-fold product and $c(m,l)$  are the coefficients of 
the elliptic genus $\sum_{m,l}c(m,l)y^lq^m$ of $X$. 
The problem was that the orbifold elliptic genus was defined in 
physical terms, and the arguments given in \cite{DMVV} did not lend itself 
to a translation into a mathematical proof.
\par The two variable elliptic genus is a very compelling invariant
for the discussion of which we refer to \cite{borlibg}.
Here we just note that it is a holomorphic function on the product 
of $\bf C$ and the upper 
half plane, which is attached to an (almost) complex manifold and is a 
weak Jacobi form if the manifold is Calabi-Yau. 
For Calabi-Yau manifolds of a dimension smaller than 12 or equal to 13,
the elliptic genus can be expressed in terms of Hirzebruch $\chi_y$ genus,
but in general, the former contains more information than the latter.
In all dimensions  elliptic genus specializes into Hirzebruch 
$\chi_y$ genus and in particular into topological Euler characteristic,
holomorphic Euler characteristic, signature, etc.
Special cases of the formula \eqref{dmvv} for these invariants 
have been proved mathematically for some time. For example,
it was shown in \cite{hirzhofer},
using the Macdonald formula \cite{Macdonald}, that if a finite group 
$G$ acts on a manifold $X$ and 
\begin{equation}
e_{orb}(X,G):={1 \over {\vert G \vert}}\sum_{fg=gf} e(X^{f,g}),
\label{orbifoldeuler}
\end{equation}
(summation is over all pairs of commuting elements; $X^{f,g}$
is the set of fixed points of both $f$ and $g$)
then:
\begin{equation}
\sum_{n=0}^{n=\infty} e_{orb}(X^n,\Sigma_n)=\prod_i {1 \over (1-t^i)^{e(X)}}.  
\label{serieseuler}
\end{equation}
On the other hand, in \cite{gottsoerg} (cf. also \cite{gottell})
it was found that the   
generating series for the $\chi_y$ genera of Hilbert schemes
of a surface $X$ is given by:
\begin{equation}
\sum_{n=0}^{n=\infty} \chi_{-y}(X^{[n]})p^n=
\exp(\sum_{m=1}^{\infty}{{\chi_{{-y}^m}(X)} \over {(1-(yp)^m)}}\,
{{p^m} \over m}).
\label{gottsche}
\end{equation}
It was observed in \cite{hirzhofer}, that in the cases when a 
crepant resolution for $X/G$ does exist, the McKay correspondence 
(cf. \cite{McKay})
can be used to prove that the Euler characteristic of 
such resolution coincides with the orbifold
Euler characteristic. This idea was used in a  more general case of 
 $\chi_y$ genus, with appropriately defined orbifold $\chi_y$ genus 
 in \cite{Batyrev.Dais}.
In the case when $X$ is a surface, the Hilbert scheme provides such a 
resolution (cf. \cite{Fogarty}) and hence the left hand side of 
(\ref{gottsche}) coincides with the generating function for
the orbifold $\chi_y$ genus of symmetric products of $X$.
Therefore, \eqref{gottsche} can be viewed as a specialization of \eqref{dmvv}.
\par This brings in the basic question: how the orbifold Euler characteristic
and  the orbifold $\chi_y$ genus, or more generally,  the orbifold elliptic genus 
of an action on a variety 
are related to the corresponding invariants of arbitrary, 
not necessarily crepant, resolution of the singularities of the orbifold. 
This question was addressed in several papers, see for example 
\cite{Batyrev}, \cite{Batyrev.Dais}, \cite{denefloeser}. 
The paper \cite{Batyrev} 
contains mathematical definitions of the orbifold $E$-function
as well as an $E$-function of singular varieties calculated via
resolutions, which is called a stringy $E$-function. 
The $E$-function of a smooth manifold is equivalent to the data given by
 Hodge numbers of the manifold, and 
it specializes to the $\chi_y$ genus. Stringy
$E$-function is defined for singular varieties with log-terminal
singularities and more generally for log-terminal pairs. 
Works \cite{Batyrev} and \cite{denefloeser} show that the 
orbifold $E$-function 
for a pair $(X,G)$ coincides with the stringy $E$-function for the pair 
$(X/G,\, {\rm image~of~ramification~divisor})$. The published version 
of \cite{Batyrev} has a gap in its canonical abelianization
algorithm, but it is now corrected by Batyrev \cite{bat.priv}.

In this paper, two notions of elliptic genus for singular varieties
are proposed. The first notion is called {\em singular elliptic genus}
and is defined for pairs (variety, divisor). Singular elliptic
genus specializes to the $\chi_y$ genus derived from the stringy
$E$-function of \cite{Batyrev}. The second notion of elliptic genus,
called {\em orbifold elliptic genus}, is defined  for
any  pair $(X,G)$ of a manifold and a finite group of its automorphisms.
Orbifold elliptic genus specializes to the $\chi_y$ genus derived 
from the orbifold $E$-function. We conjecture that the two
elliptic genera coincide for $(X/G,\,{\rm image~of~ramification~divisor})$ 
and $(X,G)$, up to an explicit normalization factor.
The advantage of the orbifold elliptic genus is that it is well-suited 
for the mathematical proof of the formula \eqref{dmvv}. On the other
hand, singular elliptic genus provides an interesting new invariant
of singular varieties. As opposed to the non-archimedian integrals over
spaces of arcs techniques of \cite{Batyrev} and \cite{denefloeser},
we use the recent result in factorization of birational maps
into a sequence of smooth blowups and blowdowns, see \cite{AKMW}.

The content of the paper is as follows. In Section \ref{basic.section},
we collect some standard definitions and results that are relevant to
the subject but may not be familiar to the reader. In Section 
\ref{sing.section} we define singular elliptic genus of a $\QQ$-Gorenstein
complex projective variety $Z$ as follows. If $\pi: Y \rightarrow Z$
is a resolution of singularities of $Z$ and 
$\alpha_k \in {\bf Q}$  are defined from the relation
$K_{Y}=\pi^*K_{Z}+\sum \alpha_kE_k$, then:
$$
\widehat{Ell}_Y(Z;z,\tau):=
\int_Y 
\Bigl(\prod_l 
\frac{(\frac{y_l}{2\pi\ii})\theta(\frac{y_l}{2\pi\ii}-z)\theta'(0)} 
{\theta(-z)\theta(\frac{y_l}{2\pi\ii})} 
\Bigr)\times 
\Bigl(\prod_k 
\frac{\theta(\frac{e_k}{2\pi\ii}-(\alpha_k+1)z)\theta(-z)} 
{\theta(\frac{e_k}{2\pi\ii}-z)\theta(-(\alpha_k+1)z)} 
\Bigr) 
$$ 
where $\theta(z,\tau)$ is the Jacobi theta function,
$y_l$  are Chern roots of $Y$ and $e_k=c_1(E_k)$.
It is shown that $\widehat{Ell}_Y(Z;z,\tau)$ 
depends only on $Z$ (rather than on the desingularization $Y$). 
Moreover, this definition is extended to pairs (variety, divisor),
and singular elliptic genus has transformation
properties of a Jacobi form if the pair satisfies a natural
Calabi-Yau condition.
Some difficulties arise only when some $\alpha_k$ equal $(-1)$,
but we generally do not need the log-terminality condition.
One application of singular elliptic genus is to the 
problem raised by M.Goreski and R.McPherson (cf. \cite{goreski}).
They were asking to determine which Chern numbers can 
be defined for singular spaces so that they are invariant 
under small resolutions. B.Totaro found a remarkable connection 
between this problem and the elliptic genus. In  \cite{totaro}
he shows that such Chern numbers must be among the combinations
of the coefficients of the two variable elliptic genus, by showing that 
these are the only Chern numbers invariant under the classical flops.
As a corollary of our definition of singular elliptic genus, we show 
that elliptic genera of any two ${\rm IH}$-small resolutions 
(or more generally two crepant resolutions) of 
a singular variety coincide, which in a sense completes the 
paper of Totaro. Unfortunately, most varieties do
not admit such resolutions, and it appears that 
Chern numbers may not be a good invariant to look for, because 
singular elliptic genera generally do not lie in the span
of the elliptic genera of smooth varieties. However, coefficients
of Taylor expansions of elliptic genera do provide an analog 
of Chern numbers for singular varieties. 

In Section \ref{orb.sect} we propose a definition of an orbifold 
elliptic genus which does not use the resolution of singularities, but uses
only information about the manifold and the fixed point sets. 
Let $G$ be a finite group acting on a manifold $X$. For $h \in G$,
let $X^h$ be a connected component of the fixed point set of $h$ 
and $TX \vert_{X^h}=\oplus V_{\lambda}, \lambda \in {\bf Q} \cap [0,1)$
be decomposition into direct sum, such that $h$ acts on $V_{\lambda}$ as the 
multiplication by $e^{2 \pi \ii \lambda}$. Let $F(h,X^h \subset X)
=\sum_{\lambda} \lambda(h)$ be the fermionic shift (cf. \cite{Batyrev.Dais}, 
\cite{Zaslow}) and:
$$V_{h,X^h\subseteq X}:= \otimes_{k\geq 1} 
\Bigr[
(\Lambda^\bullet V_0^*yq^{k-1})\otimes
(\Lambda^\bullet V_0  y^{-1}q^{k})\otimes
(Sym^\bullet V_0^*q^{k})\otimes
(Sym^\bullet V_0 q^{k})\otimes
$$
$$\otimes
\bigl[
\otimes_{\lambda\neq 0} 
(\Lambda^\bullet V_\lambda^*yq^{k-1+\lambda(h)})\otimes
(\Lambda^\bullet V_\lambda  y^{-1}q^{k-\lambda(h)})\otimes
(Sym^\bullet V_\lambda^*q^{k-1+\lambda(h)})\otimes
(Sym^\bullet V_\lambda q^{k-\lambda(h)})
\bigr]
\Bigl].
$$
Then we define (cf. Section \ref{orb.sect}):
$$Ell_{orb} (X,G;y,q) := 
y^{-\dim X /2} \sum_{\{h\},X^h}   y^{F(h,X^h\subseteq X)} \frac 1 {|C(h)|}
\sum_{g\in C(h)} L(g, V_{h,X^h\subseteq X})
$$
where $\{ h \}$ is a conjugacy class in $G$, $C(h)$ is the 
centralizer of $h$ and 
$L(g,V_{h,X^h\subseteq X})=\sum_i (-1)^i \tr(g,H^i(V_{h,X^h\subseteq X}))$
is the holomorphic Lefschetz number.  
Using the Atiyah-Singer holomorphic Lefschetz theorem, orbifold 
elliptic genus can be rewritten as follows. For a pair $g,h \in G$ of 
commuting elements, let $X^{g,h}$ be a connected component 
of the set of points in $X$ fixed by both $g$ and $h$, 
$x_{\lambda}$ be the Chern roots 
of a subbundle $V_{\lambda}$ of $TX \vert_{X^{g,h}}$ on which 
both $g$ and $h$ act via the multiplication by 
$\exp(2 \pi \ii {\lambda(g)})$ and  $\exp(2 \pi \ii {\lambda(h)})$ respectively,
and let: 
$$
\Phi(g,h,\lambda,z,\tau, x):=
\frac{{\theta(\frac{x} {2 \pi \ii }+\lambda(g)-\tau \lambda(h)-z)}
}{
{\theta(\frac{x} {2 \pi \ii }+\lambda(g)-\tau \lambda(h))}}
e^{2\pi \ii z \lambda(h)z}.
$$
Then: 
\begin{equation}
E_{orb}(X,G;z,\tau)={1 \over {\vert G \vert}} 
\sum_{gh=hg} 
\Bigl(
\prod_{\lambda(g)=\lambda(h)=0}
x_{\lambda}
\Bigr)
\prod_{\lambda} \Phi(g,h,\lambda,z,\tau,x_\lambda)[X^{g,h}]
\label{commuting}
\end{equation}
This formula generalizes \eqref{orbifoldeuler} (as we mentioned already,
the latter has as a consequence \eqref{serieseuler}, as was shown in \cite{hirzhofer}).
For a thus defined orbifold elliptic genus we prove the formula of Dijkgraaf,
Moore, Verlinde and Verlinde \eqref{dmvv}. We also show that if 
$X$ is a Calabi-Yau manifold, then $E_{orb}(X,G;z,\tau)$ is 
a weak Jacobi form.

In Section \ref{comp.gen} we conjecture (see \ref{mainconj})
that the two notions of
elliptic genera coincide, which would extend the results of 
\cite{Batyrev} and \cite{denefloeser}. We prove this conjecture
for the toric case and in dimension one. 
For Calabi-Yau hypersurfaces in Gorenstein toric Fano varieties 
the elliptic genus was defined already in \cite{borlibg}, using the 
description of the cohomology of chiral de Rham complex ${\cal MSV}$
for such hypersurfaces from \cite{Bvert} and borrowing  
the definition of elliptic genus via chiral de Rham complex in 
the nonsingular case:
$$
Ell(X)=y^{\dim X/2}{\rm Supertrace}_{H^*({\cal MSV(X)})}y^{J[0]}q^{L[0]}.
$$
Here ${\cal MSV}$ is the chiral de Rham complex constructed in 
\cite{MSV} and $J[0]$ and $L[0]$ are the operators of the  
$N=2$ super-Virasoro algebra acting  on $H^*({\cal MSV(X)})$.
We use the combinatorial description of this genus, proved in
\cite{borlibg}, and the calculation of \cite{BorGun} to show
that it coincides with the singular elliptic genus, up to 
an explicit normalization factor.

We continue to discuss Conjecture \ref{mainconj} in Section 
\ref{cobordism.section}. We show that both notions 
of elliptic genera are invariant under complex cobordisms of
$G$ action. By using the known result about cobordism classes of
the action of a cyclic group of prime order $p$, 
we prove Conjecture \ref{mainconj} for involutions.

The authors wish to thank Burt Totaro for his helpful comments.

\section{Preliminaries}\label{basic.section}

\subsection{Elliptic genus} Let $X$ be a compact (almost complex) manifold.
For a bundle $V$ on $X$ we consider the following elements in the
ring of formal power series over $K(X)$: 
$$
S_t(V)=\sum_i S^i(V)t^i \ \ \
\Lambda_t(V)=  \sum_i \Lambda^i(V)t^i
$$
where $S^i$ (resp. $\Lambda^i$) is the $i$-th symmetric (resp.
 exterior) power of $V$. 
\par \noindent Let $T_X$ (resp. $\bar T_X$) be tangent (resp. 
cotangent) bundle. The elliptic genus of $X$ can be defined as:
$$Ell(X;y,q)=\int_X ch(\Ell_{y,q})td(X)$$ where
$$\Ell_{y,q}:=y^{-{d \over 2}} \otimes_{n \ge 1} 
\Bigl(
\Lambda_{-yq^{n-1}}\bar 
T_X \otimes \Lambda_{-y^{-1}q^n} T_{X} \otimes S_{q^n}\bar T_X \otimes
S_{q^n} T_X
\Bigr).
$$
If $x_i$ are the Chern roots of $X$, i.e. for the total Chern class
we have $c(X)=\prod_i (1+x_i)$, then
\begin{equation}
Ell(X;y,q)=\int_X   \prod_i x_i {{\theta ({{x_i} \over {2 \pi \ii}}-z,\tau)} \over
{\theta ({{x_i} \over {2 \pi \ii }}, \tau)}}
\label{elltheta}
\end{equation}
where  $q=\ee^{2 \pi \ii  \tau}$ and $y=\ee^{2 \pi \ii  z}$.
In (\ref{elltheta}) 
$$
\theta(z,\tau)=q^{1 \over 8}  (2 \sin \pi z)
\prod_{l=1}^{l=\infty}(1-q^l)
 \prod_{l=1}^{l=\infty}(1-q^l \ee^{2 \pi \ii z})(1-q^l \ee^{-2 \pi \ii
z})
$$
is the Jacobi theta-function (\cite{Chandra}).
\par For $q=0$ we have: $Ell(X;y,q=0)=y^{-\frac d2} \chi_{-y}(X)$
where 
$$\chi_y(X)=\sum_{p,q} (-1)^q\dim H^q(X, \Omega^p_X)y^p$$ 
is  Hirzebruch $\chi_y$-genus (cf. \cite{Hirz1}).
In particular, $Ell(X;y=1,q=0)$ is the topological Euler characteristic,
$(-1)^{d/2}Ell(X;y=-1,q=0)$ is the signature, etc.

If $X$ is a Calabi Yau, i.e. $K_X\sim 0$, then $Ell(X;y,q)$ 
is a weak Jacobi form. Recall (cf. \cite{EZ}, \cite{Gritsenko}) 
that weight $k \in {\bf Z}$ and index 
$r \in {{1 \over 2} \bf Z}$  weak Jacobi form is a function on 
$H \times {\bf C}$ that satisfies:
$$
\phi({{a\tau +b} \over {c \tau +d}},{z \over {c \tau +d}})=
 (c \tau +d)^k e^{2 \pi i{{r cz^2} \over {c \tau +d}}  }\phi (\tau, z)
$$
$$
\phi(\tau, z+m\tau+n)=
(-1)^{2r(\lambda+\mu)} e^{-2 \pi i r (m^2 \tau +2 m z)} \phi (\tau, z)
$$
and has a Fourier expansion $\sum_{l,m}c_{m,l}y^lq^m$ with 
nonnegative $m$.

\subsection{Log-terminal singularities} We recall basic definitions
related to singular varieties.
Let $Z$ be a normal irreducible projective variety. $\QQ$-Weil
(resp. $\QQ$-Cartier) 
divisor is a linear combination with rational coefficients 
of codimension one subvarieties
(resp. Cartier divisors) on $Z$.

The canonical divisor $K_{Z}$  of $Z$
is a Weil divisor $div(s)$  where 
$s=df_1 \wedge ... \wedge df_{dim Z}$ ($f_i$ are meromorphic
functions) is a non zero rational 
differential on $Z$. We call $Z$  Gorenstein (resp. $\QQ$-Gorenstein)
is $K_Z$ is Cartier (resp. $\QQ$-Cartier).

A resolution of singularities of a variety $Z$ 
is a proper birational morphism
$f: Y \rightarrow Z$ where $Y$ is smooth. 

\begin{definition} (cf. \cite{goreski}) An  ${\rm IH}$-small resolution
of $Z$ is a regular map $Y \rightarrow Z$ such that 
for every $i \ge 1$ the set of points $z \in Z$ such that 
$\dim(f^{-1}(z)) \ge i$ has codimension greater than $2i$ in $Z$.
\end{definition}

\begin{definition} $Z$ has at worst log-terminal singularities 
if the following two conditions hold.
\par \noindent (i) $Z$ is $\QQ$-Gorenstein.
\par \noindent (ii) For a resolution $f: X \rightarrow Z$, whose
exceptional set is a divisor with simple normal crossings,
in the relation 
$K_X=f^*K_Z+\sum \alpha_iE_i$ one has $\alpha_i >-1$.
\end{definition}
A well-known result of birational geometry, see for example
\cite{KMM}, states that for any resolution of a log-terminal
variety $Z$, the coefficients $\alpha_i$ (called {\em discrepancies})
are bigger than $(-1)$. Similar definition of log-terminality exists 
for pairs $(Z,D)$ where $D$ is a $\QQ$-Weil divisor on a normal
variety $Z$ such that $(K_Z+D)$ is $\QQ$-Cartier.

\subsection{G-bundles}

Let $X$ be a complex manifold and $G$ a finite group of 
holomorphic transformations acting on $X$. Let $V$ be 
a holomorphic G-bundle on $X$, i.e. the action of $G$ on $X$
is extended to the action on $V$.
The holomorphic Lefschetz number of $g \in G$ is 
$$
L(g,V)=\sum_i (-1)^i \tr(g, H^i(X,V))
$$
Let $V^G$ be the sheaf which sections over open sets are the 
$G$-invariants of the sections of $V$.
We have (cf. \cite{Grothendieck}, spectral sequences degenerate
due to finiteness of $G$):
$$
\chi(V^G)={1 \over {\vert G \vert}} \sum_{g \in G} L(g,V).
$$
The Lefschetz numbers are given by the data around
the fixed point sets (cf. \cite{AtiyahSinger}) as follows.
Let $N^g$ be the normal bundle to the fixed point
set $X^g$ of $g$, and let ${N^g}^*$ be its dual.
In the case when the action of $G$ on a space $Y$ is 
trivial, we have $K_G(Y)=K(Y) \otimes R(G)$ (cf. \cite{AtiyahSinger})
and hence one can define $W(g) \in K(Y)$
corresponding  to 
$W \in K_G(Y)$.
In these notations:
\begin{equation}L(g,V)=
{{ch V \vert_{X^g}(g) td(T_{X^g}) } \over {ch \Lambda_{-1} 
 (N^g)^*(g)}}[X^g].
\label{holomorphicL}
\end{equation}
For $g \in G$ the normal bundle $N_{X^g}$ 
to the fixed point set $X^g$ can be decomposed into direct sum 
$N_{X^g}=\oplus_i N(\theta_i),\, \theta_i \in \QQ$ 
where each $N(\theta_i)$
is the subbundle on which $g$ acts as multiplication
by $e^{2 \pi i \theta_i}$. 
If $x_{\theta_i,j}$ are the 
Chern roots of $N(\theta_i)$ i.e. $c(N(\theta_i))=\prod_j
(1+x_{\theta_i,j})$ then (\ref{holomorphicL}) can be rewritten as: 
$$
L(g,V)={{ch (V \vert_{X^g})} \over {\prod_{i,j}
(1-e^{-x_j-\theta_{i,j}})}}td(X^g)[X^g].
$$

\section{Singular elliptic genus}\label{sing.section}
In this section we define {\em singular elliptic genus} for a large class
of singular varieties and more generally for pairs consisting of a variety and
a $\QQ$-Cartier divisor on it. This is by far the most 
general definition of elliptic genus for singular varieties constructed
to date. All varieties are assumed to be proper over ${\rm Spec}(\CC)$.

\begin{definition}\label{defellgen}
{
Let $Z$ be  a $\QQ$-Gorenstein variety 
with log-terminal singularities, and let $\pi:Y\to Z$ be 
a desingularization of $Z$ whose exceptional divisor
$E=\sum_kE_k$ has simple normal crosings.
The discrepancies $\alpha_k$ of the components $E_k$
are determined by the formula 
$$K_Y=\pi^*K_Z+\sum_k\alpha_kE_k.$$
We introduce Chern roots $y_l$ of $Y$ by $c(TY)=\prod_l(1+y_l)$
and define cohomology classes $e_k:=c_1(E_k)$. 
{\em  Singular elliptic genus} of $Z$ is then defined as 
a function of two variables $z,\tau$ given by
$$
\widehat{Ell}_Y(Z;z,\tau):=
\int_Y
\Bigl(\prod_l
\frac{(\frac{y_l}{2\pi\ii})\theta(\frac{y_l}{2\pi\ii}-z)\theta'(0)}
{\theta(-z)\theta(\frac{y_l}{2\pi\ii})}
\Bigr)\times
\Bigl(\prod_k
\frac{\theta(\frac{e_k}{2\pi\ii}-(\alpha_k+1)z)\theta(-z)}
{\theta(\frac{e_k}{2\pi\ii}-z)\theta(-(\alpha_k+1)z)}
\Bigr)
$$
where $\theta(z,\tau)$ is the Jacobi theta function, see
\cite{Chandra}. We will often suppress the $\tau$-dependence 
in our formulas.
}
\end{definition}
We will usually abuse notation and consider $\widehat{Ell}$
to be a function of 
$y=\ee^{2\pi\ii z}$ and  $q=\ee^{2\pi\ii \tau}$.
Strictly speaking, this function will be multi-valued, because
rational powers of $y$ may occur.

The key result of this section is the following theorem.
\begin{theorem}\label{main}
{
The above defined $\widehat{Ell}_Y(Z;y,q)$ does not depend on
the choice of desingularization $Y$ and therefore defines an invariant
of $Z$, which we simply denote by $\widehat{Ell}(Z;y,q)$.
}
\end{theorem}

\begin{proof} Because of the Weak Factorization Theorem of \cite{AKMW}
it suffices to show that $\widehat{Ell}_Y(Z;y,q)=\widehat{Ell}_{\tilde Y}
(Z;y,q)$ when $\tilde Y$ is obtained from $Y$ by a blowup along a 
nonsingular subvariety $X$.  We remark that the algorithm of \cite{AKMW}
is compatible with the normal crossing condition
(cf. \cite{AKMW}, Theorem 0.3.1), so we may assume that $X$ has normal
crossings with the components of the exceptional divisor of $\pi:Y\to Z$.

We will use the notations of Fulton \cite{Fulton}
for the blowup diagram 
$$
\begin{array}{ccc}
                 &   j             &             \\
    ~ \tilde X   & \longrightarrow &\tilde Y    \\
     g\downarrow~&                 &~\downarrow f \\
            ~X   & \longrightarrow &     Y     \\
                 &   i             &          
\end{array}
$$
where $\tilde X$ is the exceptional divisor of the blowdown morphism.
We also have $\pi:Y\to Z$ and $\pi\circ f:\tilde Y\to Z$. Discrepancies of
the exceptional divisors of these morphisms are related by
$$K_Y=\pi^*K_Z+\sum_k\alpha_kE_k$$
$$K_{\tilde Y}=f^*\pi^*K_Z+\sum_{k}\alpha_kE_k'+(\sum_{k}\alpha_k\beta_k
+r-1)\tilde X
$$
where $\beta_k$ is the multiplicity of $E_k$ along $X$ and 
$r$ is the codimension of $X$ in $Y$.

We will use for a while the following technical assumption.
\begin{equation}\label{technical}
\begin{array}{c}
{\rm The~normal~bundle~} N {\rm ~to~} X 
{\rm  ~inside~} Y {\rm ~ is~a~pullback~under~} i\\
{\rm ~of~some~rank ~}r{\rm ~ bundle~} M {\rm~ on~} Y.
\end{array}
\end{equation}

We have the following exact sequences of coherent sheaves on $\tilde Y$, 
see Section 15.4 of \cite{Fulton}.
$$
0\to T\tilde Y\to f^*TY\to j_*F\to 0
$$
$$
0\to j_*{\cal{O}}_{\tilde X}(-1)\to j_*g^*i^* M \to j_* F\to 0
$$
$$
0\to {\cal{O}} \to {\cal{O}}(\tilde X)\to j_*{\cal{O}}_{\tilde X}(-1)\to 0 
$$
$$
0\to f^*M(-\tilde X)\to f^*(M)\to j_*g^*i^*M\to 0
$$
Here $F$ is the tautological quotient bundle on $\tilde X$. This implies 
$$
c(T\tilde Y)=c(f^*TY)\cdot(1+\tilde x)\cdot \prod_l\frac{(1+f^*m_l-\tilde x)}
{(1+f^*m_l)}
$$
where $c(M)=\prod_l(1+m_i)$ and $\tilde x=c_1({\cal{O}}(\tilde X))$. Note also
that $c_1(E_k')=f^*e_k-\beta_k\tilde x$. Therefore,
$$
\widehat{Ell}_{\tilde Y}(Z;y,q) =  \int_{\tilde Y} 
\Bigl(\prod_l
\frac{(\frac{f^*y_l}{2\pi\ii})\theta(\frac{f^*y_l}{2\pi\ii}-z)\theta'(0)}
{\theta(-z)\theta(\frac{f^*y_l}{2\pi\ii})}
\Bigr)\times
\Bigl(
\frac{(\frac{\tilde x}{2\pi\ii})\theta(\frac{\tilde x}{2\pi\ii}-z)\theta'(0)}
{\theta(-z)\theta(\frac{\tilde x}{2\pi\ii})}
\Bigr)\times
$$
$$
\times
\Bigl(\prod_l
\frac
{\theta(\frac{f^*m_l-\tilde x}{2\pi\ii}-z)(\frac{f^*m_l-\tilde x}{2\pi\ii})
\theta(\frac {f^*m_l}{2\pi\ii})}
{\theta(\frac{f^*m_l-\tilde x}{2\pi\ii})(\frac{f^*m_l}{2\pi\ii})
\theta(\frac {f^*m_l}{2\pi\ii}-z)}
\Bigr)\times
\Bigl(\prod_k
\frac{\theta(\frac{f^*e_k-\beta_k\tilde x}{2\pi\ii}-(\alpha_k+1)z)\theta(-z)}
{\theta(\frac{f^*e_k-\beta_k\tilde x}{2\pi\ii}-z)\theta(-(\alpha_k+1)z)}
\Bigr)\times
$$
$$
\times
\Bigl(
\frac{\theta(\frac{\tilde x}{2\pi\ii}-(\alpha_{\tilde X}+1)z)\theta(-z)}
{\theta(\frac{\tilde x}{2\pi\ii}-z)\theta(-(\alpha_{\tilde X}+1)z)}
\Bigr)
$$
where $\alpha_{\tilde x} = r-1 + \sum_k \alpha_k\beta_k$.

We will now use $\int_{\tilde Y}a=\int_Yf_*(a)$. We write the Taylor expansion 
$\sum_n R_n(y,q)\tilde x^n$ of the expression under $\int_{\tilde Y}$ in
the above identity. Observe that $f_*R_0(y,q)$ is exactly the 
class in $A(Y)$ whose integral is $\widehat{Ell}_Y(Z;y,q)$, so we need to show 
that the contribution of the  rest of the terms vanishes. Notice that 
$f_*\tilde x^n=0$ for $1\leq n\leq r-1$ and $f_*\tilde x^{r+n}=i_*(s_k(i^*M))(-1)^{n+r-1}$
where $\sum_{n\geq 0}s_n t^n$ is Segre polynomial of a vector bundle, see \cite{Fulton}.
Hence, one needs to calculate 
$$
\hspace{37pt}
\int_{Y}\sum_{n\geq 0} i_*s_n(i^*M)(-1)^{n+r-1}\cdot({\rm Coeff.~at~}t^{r+n})
\Bigl[\Bigr.
\Bigl(\prod_l
\frac{(\frac{y_l}{2\pi\ii})\theta(\frac{y_l}{2\pi\ii}-z)\theta'(0)}
{\theta(-z)\theta(\frac{y_l}{2\pi\ii})}
\Bigr)\times
$$
\begin{equation}\label{twostar}
\times
\Bigl(
\frac{(\frac{t}{2\pi\ii})\theta(\frac{t}{2\pi\ii}-(\alpha_{\tilde X}+1)z)\theta'(0)}
{\theta(\frac{t}{2\pi\ii})\theta(-(\alpha_{\tilde X}+1)z)}
\Bigr)\times
\Bigl(\prod_l
\frac
{\theta(\frac{m_l-t}{2\pi\ii}-z)(\frac{m_l-t}{2\pi\ii})
\theta(\frac {m_l}{2\pi\ii})}
{\theta(\frac{m_l-t}{2\pi\ii})(\frac{m_l}{2\pi\ii})
\theta(\frac {m_l}{2\pi\ii}-z)}
\Bigr)\times
\end{equation}
$$
\times
\Bigl(\prod_k
\frac{\theta(\frac{e_k-\beta_k t}{2\pi\ii}-(\alpha_k+1)z)\theta(-z)}
{\theta(\frac{e_k-\beta_k t}{2\pi\ii}-z)\theta(-(\alpha_k+1)z)}
\Bigr)
\Bigl.\Bigr].\hfill\hspace{120pt}
$$
We denote $n_l=i^*m_l$, $f_k=i^*e_k$ and use the fact that 
$$\sum_{n\geq 0} s_n(i^*M)(-1)^nt^{-k}=\frac{t^r}{\prod_l(t-n_l)}$$
to rewrite \eqref{twostar} as 
$$
{\rm const.}\int_X ({\rm Coeff.~at~}t^{-1})
\Bigl[\Bigr.
\Bigl(\prod_l
\frac{(\frac {x_l}{2\pi\ii})\theta(\frac{x_l}{2\pi\ii}-z)\theta'(0)}
{\theta(-z)\theta(\frac{x_l}{2\pi\ii})}
\Bigr)\times
\hspace{40pt}
$$
\begin{equation}\label{threestar}
\times
\Bigl(
\frac{\theta(\frac{t}{2\pi\ii}-(\alpha_{\tilde X}+1)z)\theta'(0)}
{\theta(\frac{t}{2\pi\ii})\theta(-(\alpha_{\tilde X}+1)z)}
\Bigr)\times
\Bigl(\prod_l
\frac
{\theta(\frac{n_l-t}{2\pi\ii}-z)
\theta(\frac {n_l}{2\pi\ii})}
{\theta(\frac{n_l-t}{2\pi\ii})(\frac{n_l}{2\pi\ii})
\theta(\frac {n_l}{2\pi\ii}-z)}
\Bigr)\times
\end{equation}
$$
\times
\Bigl(\prod_k
\frac{\theta(\frac{f_k-\beta_k t}{2\pi\ii}-(\alpha_k+1)z)\theta(-z)}
{\theta(\frac{f_k-\beta_k t}{2\pi\ii}-z)\theta(-(\alpha_k+1)z)}
\Bigr)
\Bigl.\Bigr].
\hspace{110pt}
$$
Here we denote $c(TX)=\prod_l(1+x_l)$ and use $c(TX)=i^*c(TY)/i^*c(M)$.
To show that \eqref{threestar} 
is zero, observe that the function whose coefficient 
at $t^{-1}$ is measured, is elliptic in $t$. Really, 
$t\to t+2\pi\ii$ obviously keeps it unchanged,
and $t\to t+ 2\pi\ii\tau$ does not change it, because of $\alpha_{\tilde X}=
\sum_k\alpha_k\beta_k +r-1$. We have used here the fact that none of 
the $\alpha$-s is equal to $(-1)$, which follows from the condition that 
$Z$ is log-terminal, see for instance \cite{KMM}. It remains to show 
that $t=0$ is the only pole
of this function up to the lattice $2\pi\ii(\ZZ+\ZZ\tau)$, so the residue 
is zero. To do so, observe that the normal crossing condition implies 
$\beta_k\in\{0,1\}$, and moreover, whenever $\beta_k=1$ the corresponding 
factor 
$
\theta(\frac{f_k-t}{2\pi\ii}-z)
$
in the denominator of the last product is offset by a factor 
$
\theta(\frac{n_l-t}{2\pi\ii}-z)
$
in the numerator of the second product.

We will now get rid of the assumption \eqref{technical}. 
Indeed, it is easy to see that
the difference between $\widehat{Ell}_Y(Z;y,q)$ and 
$\widehat{Ell}_{\tilde Y}(Z;y,q)$ can be written 
as a degree of an element of $A(\tilde X)$ that is preserved when one
deforms $i:X\to Y$ to the embedding of $X$ into the normal cone, for
which the assumption \eqref{technical} is satisfied.
\end{proof}

\begin{remark}
{
We have not significantly used the log-terminality condition, except for
the fact that we did not have to divide by $\theta(0\cdot z)$. Therefore,
singular elliptic genera can in fact be defined for all varieties
that admit a resolution with no $(-1)$ discrepancies. In fact, we 
will extend our definition to the category of pairs that consist
of an algebraic variety and a $\QQ$-Cartier divisor on it.
}
\end{remark}

\begin{definition}\label{singelllog}
{ 
Let $Z$ be a projective variety, and let $D$ be an arbitrary $\QQ$-Weil
divisor such that $K_Z+D$ is a $\QQ$-Cartier divisor on 
$Z$. Let $\pi:Y\to Z$ be a desingularization of $Z$.
We denote by $E=\sum_kE_k$ the exceptional divisor of $\pi$ plus
the sum of the proper preimages of the components of $D$ and assume that
it has simple normal crossings. The 
discrepancies $\alpha_k$ of the components $E_k$
are determined by the formula 
$$K_Y=\pi^*(K_Z+D)+\sum_k\alpha_kE_k$$
and the requirement that the discrepancy of the proper transform
of a component of $D$ is the opposite of the coefficient of
$D$ at that component.
We introduce Chern roots $y_l$ of $Y$ by $c(TY)=\prod_l(1+y_l)$
and define
$$
\widehat{Ell}_Y(Z,D;y,q):=
\int_Y
\Bigl(\prod_l
\frac{(\frac{y_l}{2\pi\ii})\theta(\frac{y_l}{2\pi\ii}-z)\theta'(0)}
{\theta(-z)\theta(\frac{y_l}{2\pi\ii})}
\Bigr)\times
\Bigl(\prod_k
\frac{\theta(\frac{e_k}{2\pi\ii}-(\alpha_k+1)z)\theta(-z)}
{\theta(\frac{e_k}{2\pi\ii}-z)\theta(-(\alpha_k+1)z)}
\Bigr)
$$
where as usual $y=\ee^{2\pi\ii z}$, $q=\ee^{2\pi\ii \tau}$, the 
$\tau$-dependence is suppressed, and $e_k=c_1(\cal{O}(E_k))$.
If some of the discrepancies $\alpha_k$ equal $(-1)$, then 
we try to define the elliptic genus as follows. For any ample 
effective Cartier divisor $H$ on $Z$ that contains all
singular points of $Z$ and all components of $D$ we calculate 
$$
\lim_{n\to\infty} \widehat{Ell}_Y(Z,D+H/n;z,\tau)
$$
for each $(z,\tau)$. If such limit exists and is independent of 
$H$, then we call it $\widehat{Ell}_Y(Z,D;z,\tau)$.
Notice that if $n$ is sufficiently big, then the discrepancies of
all divisors $E_k$, calculated for the pair $(Z,D+H/n)$ 
are not equal to $(-1)$.
It is also  easy to see that if there are no $(-1)$ discrepancies,
then $\lim_{n\to\infty} \widehat{Ell}_Y(Z,D+H/n;z,\tau)=
\widehat{Ell}_Y(Z,D;z,\tau)$.}
\end{definition}

\begin{theorem}\label{mainlog}
{
The above defined elliptic genus does not depend on the choice 
of the desingularization $\pi:Y\to Z$. We will therefore denote it
simply by $\widehat{Ell}(Z,D;y,q)$.
}
\end{theorem}

\begin{proof} 
Any two resolutions of singularities of $Z$ can be connected by a 
sequence of blowups and blowdowns. Let $Y$ and $\tilde Y$ be 
two resolutions of $Z$, such that $\tilde Y$ is the blowup of $Y$ 
as in the proof of Theorem \ref{main}.  For any $H$ and $n$ big 
enough to assure that all discrepancies are not equal to $(-1)$, 
the proof
of Theorem \ref{main} implies that 
$$
\widehat{Ell}_Y(Z,D+H/n;z,\tau)=\widehat{Ell}_{\tilde Y}(Z,D+H/n;z,\tau).
$$
Then elliptic genera defined via $Y$ and $\tilde Y$ coincide by 
definition. 
\end{proof}

\begin{remark}
{
The reason behind extending the definition of the elliptic genus
via above limits is the following. In the non-log-terminal case,
it is conceivable that there exist two resolutions of singularities
without $(-1)$ discrepancies that can only be connected via
resolutions with $(-1)$ discrepancies. The above theorem assures 
that elliptic genera defined via such resolutions are the same.
}
\end{remark}

\begin{proposition}
{ 
Elliptic genera of two different crepant resolutions of a 
Gorenstein projective variety coincide.
}
\label{twocrepant}
\end{proposition}

\begin{proof} 
We will show that elliptic genus of a crepant resolution $Y$
of a variety $X$ equals to the singular elliptic genus of $X$.
If the exceptional set of the morphism $\pi:Y\to X$ 
is a divisor with simple normal crossings, then it is enough
to observe that in the Definition \ref{defellgen} the second product
is trivial. In general, we can further blow up $Y$ to get $\mu:Z\to Y$
so that the exceptional sets of $\mu$ and  $\pi\circ\mu:Z\to X$ are divisors
with simple normal crossings. Then singular elliptic genera of $Y$ 
and $X$ calculated via $Z$ are given by the same formula, because the
discrepancies coincide.
\end{proof}

\begin{remark}
{
In particular, the above proposition shows that the statement 
of Theorem 8.1 of \cite{totaro} can be extended to the full elliptic
genus.
}
\end{remark}

The following proposition shows that when $q\to 0$, 
we recover a formula for $\chi_y$ genus of $(Z,D)$ which follows from 
\cite{Batyrev}.

\begin{proposition}\label{normcross}
{ Let $(Z,D)$ be a log-terminal pair, see \cite{Batyrev}.
Then
$$
\widehat{Ell}(Z,D;u,q=0) = (u^{-\frac 12} - u^{\frac 12})^{\dim Z}
 E_{st}(Z,D;u,1)
$$
where $E_{st}$ is defined in \cite{Batyrev}.
}
\end{proposition}

\begin{proof}
To avoid confusion, we immediately remark that the second arguments 
in singular elliptic genus and in Batyrev's $E$-function have drastically 
different meanings.
The definition of $E_{st}(Z,D)$ in \cite{Batyrev} could be stated as 
$$E_{st}(Z,D;u,v):=\sum_{J\subset I} E(E_J;u,v) \prod_{j\in J}
\Bigl(\frac{uv-1}{uv^{\alpha_j+1}-1}-1\Bigr)$$
where $\sum_{i\in I} \alpha_i E_i$ is the exceptional divisor of
a resolution $Y\to Z$ together with proper preimages of the components of 
$D$, and  is assumed to have normal crossings. Polynomials $E(E_J;u,v)$
are defined in terms of mixed Hodge structure on the cohomology
of $E_J$, see \cite{Batyrev}.
Subvariety $E_J$ is $\cap_{j\in J}E_j$, and the sum includes the empty
subset $J$.

For each $J$ 
$$
E_{st}(E_J;u,1)=\int_{E_j}\prod_{i=1}^{\dim E_J} 
\frac {(1-u\ee^{-x_{i,J}})x_{i,J}}{1-\ee^{-x_{i,J}}},
$$
where $c(TE_J)=\prod_i(1+x_i,J)$. The adjunction formula
for complete intersections yields
$$c(TE_J)=i_J^*(c(TY))/\prod_{j\in J}(1+i_J^*c_1(E_j)),$$
where $i_J:E_J\to Y$ is the closed embedding. We then obtain 
$$
E(E_J;u,1)=\int_{E_j}\prod_{i=1}^{\dim Y} 
\frac {(1-u\ee^{-i_J^*x_{i}})i_J^*x_{i}}{1-\ee^{-i_J^*x_{i}}}
\prod_{j\in J} \frac {1-\ee^{-i_J^*e_j}}{(1-u\ee^{-i_J^*e_j})i_J^*e_j}
$$
$$
=\int_Y \prod_{i=1}^{\dim Y}
\frac {(1-u\ee^{-x_{i}})x_{i}}{1-\ee^{-x_{i}}}
\prod_{j\in J} \frac {1-\ee^{-e_j}}{1-u\ee^{-e_j}}
$$
where $c(TY)=\prod_i(1+x_i)$. When we plug this result into
Batyrev's formula, we get
$$
E_{st}(Z,D;u,1)=
\int_Y \prod_{i=1}^{\dim Y}
\frac {(1-u\ee^{-x_{i}})x_{i}}{1-\ee^{-x_{i}}}
\prod_{j\in I}\Bigl( 1+ \frac {(u-u^{\alpha_j+1})(1-\ee^{-e_j})}
{(u^{\alpha_j+1}-1)(1-u\ee^{-e_j})}\Bigr)
$$
$$
\hspace{37pt}
=
\int_Y \prod_{i=1}^{\dim Y}
\frac {(1-u\ee^{-x_{i}})x_{i}}{1-\ee^{-x_{i}}}
\prod_{j\in I}
\frac {(u-1)(1-u^{\alpha_j+1}\ee^{-e_j})}{(u^{\alpha_j+1}-1)(1-u\ee^{-e_j})}
$$
$$
= 
(u^{-\frac 12} - u^{\frac 12})^{-\dim Z}
\lim_{q\to 0} \widehat{Ell}(Z,D;u,q).
\hspace{20pt}
$$
\end{proof}

The following simple proposition establishes modular properties of 
singular elliptic genus in Calabi-Yau case.

\begin{proposition}
{
Let $(Z,D)$ be a Calabi-Yau pair, in the sense that 
$K_Z+D$ is zero as a $\QQ$-Cartier divisor.
If $Z$ admits a resolution such that all the discrepancies
of $(Z,D)$ are not equal to $(-1)$, then 
singular elliptic genus $\widehat{Ell}(Z,D;y,q)$ 
has transformation properties of Jacobi form of 
weight $\dim Z$ and index $0$ for the subgroup of the full
Jacobi group generated by 
$$(z,\tau)\to(z+n,\tau),~(z,\tau)\to(z+n\tau,\tau),~(z,\tau)\to(z,\tau+1),
~(z,\tau)\to(z/\tau,-1/\tau)$$
where $n$ is the least common denominator of the discrepancies.
}
\end{proposition}

\begin{proof}
Transformation properties of
$\theta(z,\tau)$ under $(z,\tau)\to(z+1,\tau)$ and 
$(z,\tau)\to(z+\tau,\tau)$ together with Calabi-Yau condition
$$K_Y=\sum_{k}\alpha_k E_k$$
assure that 
$$\widehat{Ell}(Z,D;z+n,\tau)=\widehat{Ell}(Z,D;z+n\tau,\tau)
=\widehat{Ell}(Z,D;z,\tau).$$
We needed here that $n\alpha_k\in \ZZ$.
Similarly, the transformation properties of $\theta$ under
$(z,\tau)\to(z,\tau+1)$ show that 
$$\widehat{Ell}(Z,D;z,\tau+1)=\widehat{Ell}(Z,D;z,\tau).$$
It remains to investigate what happens under
$(z,\tau)\to(z/\tau,-1/\tau)$. For this, one considers
the change $(e_k,y_l)\to (e_k/\tau,y_l/\tau)$ in
the formula of the Definition \ref{singelllog}. A rather 
lengthy but straightforward calculation, similar to that of 
Theorem 2.2 of \cite{borlibg}, shows that
$$
\widehat{Ell}(Z,D;\frac z\tau , -\frac 1\tau)=
\tau^{\dim Z}
\widehat{Ell}(Z,D;z,\tau).$$
\end{proof}

Another application of our techniques is the following theorem,
which complements similar results for Hodge numbers of Calabi-Yau
manifolds, see for example \cite{Batyrev.betti} and \cite{denefloeserinv}.
\begin{theorem}\label{birCY}
{
Elliptic genera of two birationally equivalent Calabi-Yau manifolds
coincide. Moreover, the statement is true for smooth projective
algebraic manifolds $X$ with $nK_X\sim 0$ for some $n$.
}
\end{theorem}

\begin{proof} Let $Z_1$ and $Z_2$ be two birationally equivalent Calabi-Yau
manifolds or their generalizations above. 
Let $Y$ be a desingularization of the closure of the graph
of the birational equivalence, so that $\pi_{1,2}:Y\to Z_{1,2}$ are
regular birational morphisms. Let $n$ be the smallest integer so that 
$nK_{Z_{1,2}}$ is rationally equivalent to zero, and therefore has a 
global section. Global sections of the pluricanonical bundle are birational 
invariants, so one can consider the divisor $\sum_k a_kE_k$ of this section
on $Y$. It is easy to see that for both morphisms $\pi_1$ and $\pi_2$ 
the exceptional divisor is $\sum_k (a_k/n)E_k$, which we can then assume
to have simple normal crossings (perhaps by passing to a new desingularization).
Therefore, elliptic genera of $Z_{1,2}$ are calculated on $Y$ using the 
same discrepancies.
\end{proof}

\begin{remark}
{
It is interesting to compare the results of this section with the 
work of Totaro in \cite{totaro}, who tried to see which Chern numbers
can be meaningfully defined for singular varieties. 
For varieties that admit ${\rm IH}$-small resolutions, 
singular elliptic genus provides the maximum possible 
collection of such numbers. Totaro has shown that every flop-invariant
Chern number comes from the elliptic genus and obtained
partial results in the opposite direction by means of 
intersection cohomology.

In general, coefficients of the singular elliptic genus of $Z$ at 
$y^kq^l$ provide analogs of Chern numbers of singular varieties
in the following sense.

\noindent 
1. They are the invariants of isomorphism class of singular spaces.

\noindent 
2. For manifolds these invariants are the usual Chern numbers (i.e. 
linear combinations of 
$c_{i_1}(X) \cdot \cdot \cdot c_{i_N}(X)[X]$ where
$\sum_{k=1}^{k=N} i_k =\dim X$ and $[X]$ is the fundamental class 
of a manifold $X$).

\noindent 
3. These invariants are unchanged under small resolutions.

\noindent 
In fact, for singular varieties, elliptic genera may contain
more information than in the non-singular case. 
For varieties with non-Gorenstein 
singularities, singular elliptic genus may depend on rational powers of $y$.
Moreover, there exist examples of Gorenstein varieties whose
elliptic genera do not lie in the span of elliptic genera of nonsingular
varieties. This can be seen already at the level of $\chi_y$ genus,
see an example in \cite{Batyrev.cangor} of a variety with 
Gorenstein canonical singularities whose $E$-function is not a 
polynomial. 

We hope that elliptic genera of singular varieties can be interpreted 
as nontrivial invariants of not yet defined cobordism theory of 
singular spaces. Transformations leaving the singular elliptic genus 
invariant in such theory for smooth manifolds should include usual cobordisms as well as flops. It would be interesting to compare our results 
with the invariants of Witt spaces studied by P.Siegel,
the latter however were defined  in  $SO$ rather than in complex category 
(cf. \cite{goreski}, \cite{siegel}).
}
\end{remark}

\begin{remark}
{
We do not have a good understanding of the reason why $(-1)$ discrepancies
seem to be a problem. One can observe, however, that in the case
of a surface singularity obtained by contracting a single smooth curve on
a smooth surface to a point, the discrepancy is $(-1)$ if and
only if the curve in question is elliptic.
}
\end{remark}

\section{Orbifold elliptic genus and DMVV formula}\label{orb.sect}
In this section we define {\em orbifold elliptic genus}, which 
we conjecture to equal the singular elliptic genus of 
Section \ref{sing.section}. 
We delay the comparison of two genera until Section \ref{comp.gen}.
Instead, the goal of this section is to show how this definition of 
orbifold elliptic genus allows one to recover the formula of
\cite{DMVV} whose derivation was based partly on heuristic 
string-theoretic arguments.
Our definition of elliptic genus is inspired by the calculations 
of \cite{Bvert}. 
\begin{definition}\label{orbifold.def}
{ Let $X$ be an algebraic variety acted upon by a finite group $G$.
We assume that the subgroup of elements of $G$ acting trivially on $X$ 
contains only the identity.
We define the following function of two variables that we call 
{\em orbifold elliptic genus} of $X/G$:
$$Ell_{orb} (X,G;y,q) := 
y^{-\dim X /2} \sum_{\{h\},X^h}   y^{F(h,X^h\subseteq X)} \frac 1 {|C(h)|}
\sum_{g\in C(h)} L(g, V_{h,X^h\subseteq X})
$$
where $F(h,X^h\subseteq X)$ is the fermionic shift (cf. \cite{Zaslow}, 
\cite{Batyrev.Dais})
and $V_{h,X^h\subseteq X}$ is a vector bundle over $X^h$ 
defined as follows. Let $TX|_{X^h}$ decompose into eigensheaves for 
$h$ as  
\begin{equation}
V_0\oplus(\oplus_{\lambda: <h>\to \QQ/\ZZ} V_{\lambda}).
\label{directsum}
\end{equation}
We lift $\lambda(h)$ to a rational number in $[0,1)$. Then
$V_{h,X^h\subseteq X}$
is defined as 
$$
V_{h,X^h\subseteq X}:= \otimes_{k\geq 1} 
\Bigr[
(\Lambda^\bullet V_0^*yq^{k-1})\otimes
(\Lambda^\bullet V_0  y^{-1}q^{k})\otimes
(Sym^\bullet V_0^*q^{k})\otimes
(Sym^\bullet V_0 q^{k})\otimes
$$
$$\otimes
\bigl[
\otimes_{\lambda\neq 0} 
(\Lambda^\bullet V_\lambda^*yq^{k-1+\lambda(h)})\otimes
(\Lambda^\bullet V_\lambda  y^{-1}q^{k-\lambda(h)})\otimes
(Sym^\bullet V_\lambda^*q^{k-1+\lambda(h)})\otimes
(Sym^\bullet V_\lambda q^{k-\lambda(h)})
\bigr]
\Bigl].
$$
}
\end{definition}

\begin{remark}
{
Another way to state this definition is 
$$
Ell_{orb}(X,G;y,q):=y^{-\dim X/2}\sum_{\{h\},X^h}
y^{F(h,X^h\subseteq X)}\chi(H^\bullet(V_{h,X^h\subseteq X}^{C(h)})).
$$
}
\end{remark}

\begin{theorem}
{ Let $X$ and $G$  be as above and let $X^{g,h}$ be 
the set of fixed points of a pair of commuting 
elements $g,h \in G$.
Let $TX \vert_{X^{g,h}}=\oplus W_{\lambda}$ 
be the decomposition (refinement of \eqref{directsum}) 
of the restriction on $X^{g,h}$ of
the tangent bundle into 
direct sum of line bundles on which $g$ (resp. $h$) acts
as multiplication by $\ee^{2 \pi \ii \lambda(g)}$ (resp. 
  $\ee^{2 \pi \ii \lambda(h)}$). Denote by $x_{\lambda}$ the 
Chern roots of the bundle 
$W_{\lambda}$.
\bigskip 
\par \noindent 1. We have: 
$$
Ell_{orb}(X,G)=
$$
$$
{1 \over {\vert G \vert }}
\sum_{g,h, gh=hg} 
\Bigl(
\prod_{\lambda(g)=\lambda(h)=0} x_{\lambda} 
\Bigr)
\prod_{\lambda} {{ \theta(\tau,{{x_{\lambda}} \over {2 \pi \ii }}+
 \lambda (g)-\tau \lambda(h)-z )} \over 
{ \theta(\tau,{{x_{\lambda}} \over {2 \pi \ii }}+
 \lambda (g)-\tau \lambda(h))} } e^{2 \pi \ii \lambda(h)z}
[X^{g,h}].
$$

\noindent 
2. Let $X$ be a Calabi-Yau of dimension $d$, such that $H^0(X,K_X)=\CC$.
Denote by $n$ the order of $G$ in ${\rm Aut} H^0(X,K_X)$. Then 
$Ell_{orb}(X,G)$ is a weak Jacobi form of weight $0$ and index $d/2$
with respect to subgroup of the Jacobi group $\Gamma^J$
generated by transformations
$$(z,\tau) \to (z+n,\tau),
~(z,\tau) \to (z+n\tau, \tau),
~(z,\tau)\to (z,\tau +1), 
~(z,\tau)\to (\frac z\tau,-\frac {1}\tau).$$
In particular, if the action preserves holomorphic volume then 
$Ell_{orb}(X,G)$ is a weak Jacobi form of weight $0$ and index 
$d/2$ for the full Jacobi group.
}
\end{theorem}

\begin{proof} 
We replace the contribution of each conjugacy class by an average contribution
of its elements to obtain
$$Ell_{orb}(X,G)=
\frac{1}{ {\vert G \vert}}y^{-\dim X/2}
 \sum_{gh=hg} y^{F(h,X^h \subset X)} L(g,V_{h,X^h \subset X}).$$
Using holomorphic Lefschetz theorem, we obtain: 
$$
Ell_{orb}(X,G)=\frac{1 }{ {\vert G \vert}} y^{-\dim X/2}
 \sum_{gh=hg} y^{F(h,X^h \subset X)}
\frac {{ch(V_{h,X^h \subset X} \vert_{X^{g,h}} )(g) td(T_{X^{g,h}})}[X^{g,h}]
}{ {{ch \Lambda_{-1} 
 (N^g_{X^h})^*(g)}}}, 
$$
where $N^g_{X^h}$ is a the normal bundle to $X^{g,h} $ in $X^h$.
An explicit calculation of the Chern and Todd classes then yields
$$
Ell_{orb}(X,G)=
%
%
\frac{1 }{ {\vert G \vert}} \sum_{gh=hg} y^{F(h,X^h \subset X)-\dim X/2} 
\Bigl(
\prod_{\lambda(g)=\lambda(h)=0} x_{\lambda}
\Bigr)\times
$$
$$
\times\prod_{k \ge 1,\lambda}
\frac                                                        
{(1-yq^{k-1+\lambda(h)}e^{-x_{\lambda}-2 \pi \ii \lambda(g)})
 (1-y^{-1}q^{k-\lambda(h)}e^{x_{\lambda}+2 \pi \ii \lambda(g)})} 
 {(1-q^{k-1+\lambda(h)}e^{-x_{\lambda}-2 \pi \ii \lambda(g)})
(1-q^{k-\lambda(h)}e^{x_{\lambda}+2 \pi \ii \lambda(g)})}=
$$
$$\frac{1 }{\vert G \vert} \sum_{gh=hg}
\Bigl(
\prod_{\lambda(h)=\lambda(g)=0} x_{\lambda} 
\Bigr)
\prod_{\lambda}
\frac{\theta(\frac {x_\lambda} {2 \pi \ii }+\lambda(g)-\tau \lambda(h)-z)} 
{\theta(\frac {x_\lambda} { 2 \pi \ii }+\lambda(g)-\tau \lambda(h))}
e^{2\pi \ii z\lambda(h)}[X^{g,h}]
$$
which proves the first part of the theorem.

To verify the modular property, we denote 
$$\Phi(g,h,\lambda,z,\tau, x):=
\frac{\theta(\frac {x} {2 \pi \ii }+\lambda(g)-\tau \lambda(h)-z)}
{\theta(\frac x {2 \pi \ii }+\lambda(g)-\tau \lambda(h))}
e^{2\pi \ii z \lambda(h)}$$ 
where $\lambda$ is a character of the subgroup of $G$ generated by
$g$ and $h$. Then 
\begin{equation}\label{hgphi}
E_{orb}(z,\tau)=
\frac{1}{\vert G \vert} 
\sum_{gh=hg} 
\Bigl(
\prod_{\lambda(g)=\lambda(h)=0}
x_{\lambda}
\Bigr)
\prod_{\lambda} \Phi(g,h,z,\tau,x_\lambda)[X^{g,h}] 
\end{equation}
where we suppress $(X,G)$ from the notations for the sake of brevity.
We have:
$$\Phi(g,h,\lambda,z+1,\tau,x)=-\ee^{2 \pi \ii \lambda(h)} \cdot 
\Phi(g,h,\lambda,z,\tau,x)$$
and hence $Ell_{orb}(z+n,\tau)=(-1)^{dn} Ell_{orb}(z,\tau)$,
since by assumption $n \cdot \sum \lambda(h) \in {\bf Z}$.
It is clear that
$$\Phi(g,h,\lambda,z,\tau+1,x)=\Phi(gh^{-1},h,\lambda,z,\tau,x)$$ and hence
$Ell_{orb}(z,\tau+1)=Ell_{orb}(z,\tau)$.
We have: 
$$\Phi(g,h,\lambda,z+n\tau,\tau,x)=
(-1)^n
\ee^{-2 \pi \ii nz -\pi \ii n^2 \tau}\ee^{nx+2 \pi \ii n\lambda(g)} 
\cdot 
\Phi(g,h,\lambda,z,\tau,x)$$ and hence
$$
Ell_{orb}(z+n \tau,\tau)=(-1)^{dn} e^{-2 \pi \ii dn z - \pi \ii dn^2 \tau}Ell_{orb}
(z,\tau)$$
since $X$ is Calabi-Yau and $n\lambda(g)\in\ZZ$.
Finally,
$$\Phi(g,h,\lambda,\frac z \tau,\ -\frac 1  \tau,\frac x \tau)=
{{\theta(-{z \over \tau}+{x \over {2 \pi i
 \tau} }+\lambda(g)+{{\lambda(h)} \over \tau},-\frac 1\tau)}
 \over 
{\theta({x_{\lambda} \over 2 \pi \ii }+\lambda(g)+ {\lambda(h) \over {\tau}}
,-\frac 1\tau)}}
e^{\frac {2\pi \ii z \lambda(h)}\tau}=
$$
$$
e^{{\pi \ii z^2 \over \tau} -{{2 \pi \ii z} \over {\tau}}
({x \over {2 \pi \ii }}+ \lambda(g) \tau +\lambda(h))}
{{\theta (-z +{x \over {2 \pi \ii }}+\lambda(g) \tau +{\lambda(h)},\tau)}  
 \over 
{\theta(
{x \over 2 \pi \ii }+\lambda(g) \tau +\lambda(h),\tau)}}
e^{\frac{2\pi \ii z \lambda(h)}\tau}=
$$
$$
e^{{\pi \ii z^2 \over \tau}-\frac {zx}\tau}
\cdot
{{\theta (-z +{x \over {2 \pi \ii }}+\lambda(g) \tau +{\lambda(h)},\tau)}  
 \over 
{\theta(
{x \over 2 \pi \ii }+\lambda(g) \tau + \lambda(h),\tau)}}
e^{2 \pi \ii z (-\lambda(g))}=
e^{{\pi \ii z^2 \over \tau}-\frac {zx}\tau}
\cdot\Phi(h,g^{-1},\lambda,z,\tau, x).
$$ 
Then the Jacobi transformation properties follows easily from \eqref{hgphi},
similar to the proof of Theorem 2.2 in \cite{borlibg}.

It is straightforward to see from \eqref{hgphi} that orbifold elliptic
genus is holomorphic and has the Fourier expansion with
non-negative powers of $q$.
\end{proof}

We will apply our definition of the orbifold elliptic genus 
to symmetric products of a smooth variety.
This will give a mathematical justification of the physical calculation
performed in \cite{DMVV}. More precisely, we calculate the 
generating function for the orbifold elliptic genera introduced 
above for the action of the symmetric groups. 
Our calculation to certain extent follows 
\cite{DMVV}, but we now have precise mathematical definitions.
\begin{theorem}\label{DMVVtheorem}
{
Let $X$ be a smooth variety with elliptic genus
$\sum_{m,l}c(m,l)y^lq^m$, where elliptic genus is normalized
as in \cite{DMVV} and \cite{borlibg}. Then 
$$
\sum_{n\geq 0}p^n Ell_{orb}(X^n,\Sigma_n;y,q)=
\prod_{i=1}^\infty \prod_{l,m}\frac 1
{(1-p^iy^lq^m)^{c(mi,l)}}.
$$
}
\end{theorem}

We shall start with the following lemma essentially 
contained in (\cite{DMVV}, Section 2.2), which we 
include only for completeness.
\begin{lemma}\label{DMVVlemma}  
{
Let $V=V_{\rm even} \oplus V_{\rm odd}$ be a supersymmetric space
and $A$ and $B$ be two commuting operators preserving 
parity decomposition of $V$, such that $B$ has only non-negative
integer eigenvalues. We assume that $V$ splits into
a direct sum of eigenspaces $V_{m}$ of the operator $B$,
and each $V_m$ is finite-dimensional. 
Define 
$$\chi(V)(y,q)={\rm Supertrace}_V y^Aq^B:=
\tr_{V_{\rm even}}(y^Aq^B) - \tr_{V_{\rm odd}}(y^Aq^B) =\sum_{m,l} d(m,l)q^my^l$$
where $d(m,l)$ is the superdimension of the space $V_{m,l}=\{ v \in V 
\vert Av=lv,Bv=mv \}$.
The operators $A$ and $B$ act on the space 
of invariants of the symmetric group acting on  $V^{ \otimes^N}$ 
and 
$$\sum_N p^N {\rm Supertrace}_{Sym^N(V)}y^Aq^B=
\prod_{m,l}{1 \over {(1-pq^my^l)^{d(m,l)}}}$$
where the right hand-side is expanded as a power series in $q$ and $p$.
}
\end{lemma}

\begin{proof}
It is easy to see that it is enough to check the lemma for a one-dimensional 
space $V=V_{m,l}$. If $V$ is even, then 
$$\sum_N p^N {\rm Supertrace}_{Sym^N(V)}y^Aq^B=
\sum_{N\geq 0} p^N y^{Nl}q^{Nm}= (1-pq^my^l)^{-{\rm superdim}V}.$$
If $V$ is odd, then 
$$\sum_N p^N {\rm Supertrace}_{Sym^N(V)}y^Aq^B=
1- p y^{l}q^{m} = (1-pq^my^l)^{-{\rm superdim}V}.$$ 
\end{proof}

We are  now ready to prove Theorem \ref{DMVVtheorem}.
\begin{proof}
We observe that for a fixed $k$ the conjugacy classes of $\Sigma_k$ 
are indexed by the numbers $a_i$ of cycles of length $i$ in the permutation. 
For each $h\in \Sigma_k$ the fixed point set $(X^k)^h$ consists 
of the Cartesian products of several copies of $X$, one for each cycle.
For a cycle of length $i$ the corresponding $X$ is embedded into $X^i$.
The centralizer group is a semidirect product of its normal subgroup
$\prod_{i} (\ZZ/i\ZZ)^{a_i}$ which acts by cyclic permutations inside
cycles of $h$ and the product of symmetric groups $\prod_i \Sigma_{a_i}$ 
that act by permuting cycles of the same length. 

Our definition of the elliptic genus then gives
$$
\sum_{n\geq 0}p^n Ell_{orb}(X^n,\Sigma_n;y,q)=
\sum_{a_1,a_2,\ldots,a_n} p^{a_1+2a_2+\ldots+na_n}
y^{-\frac {\dim X}2 (a_1+2a_2+\ldots+na_n)}$$
\begin{equation}\times
\prod_{i=1}^n 
y^{a_iF(i-cycle,X\subseteq X^i)}
\chi\Bigl(\bigl(H^\bullet(V_{i-cycle,X\subseteq X^i})^{\otimes a_i}
\bigr)^{\Sigma_{a_i}\rtimes (\ZZ/i\ZZ)^{a_i}}\Bigr)
\label{firststep}
\end{equation}
$$
=\prod_{i=1}^\infty 
\chi\Bigr(
Sym^\bullet\bigl(p^iy^{-i\dim X/2 + F(i-cycle,X\subseteq X^i)}
H^\bullet(V_{i-cycle,X\subseteq X^i})^{\ZZ/i\ZZ}\bigr)
\Bigr).
$$
The symbol $Sym$ here should be interpreted as the supersymmetric product
where the cohomology of $V_{h,X^h\subseteq X}$  is given parity by the 
sum of the cohomology number and the parity of the exterior algebras.

We will now calculate
$$
\chi_i(y,q)=
\chi
\bigl(p^iy^{-i\dim X/2 + F(i-cycle,X\subseteq X^i)}
H^\bullet(V_{i-cycle,X\subseteq X^i})^{\ZZ/i\ZZ}\bigr).
$$
We denote the $i$-cycle by $h$ and observe that 
$$TX^i|_X = \oplus_{j=0,\ldots,i-1;\lambda(h)=\frac ji} TX_{j}.
$$
This implies $F(h, X\subseteq X^i)= \dim X \sum_{j=0}^{i-1} \frac ji =
\dim X \frac {(i-1)}2$,
which allows us to write 
$$
\chi_i(y,q)= p^i y^{-\dim X/2}
\chi
\Bigl[
\bigl[H^\bullet(
\otimes_{k\geq 1} 
\bigr[
(\Lambda^\bullet T^*yq^{k-1})\otimes
(\Lambda^\bullet T  y^{-1}q^{k})\otimes
(Sym^\bullet T^*q^{k})
$$
$$
\otimes
(Sym^\bullet T q^{k})\otimes
\bigl[
\otimes_{j=1,\ldots,i-1} 
(\Lambda^\bullet T^*yq^{k-1+\frac ji})\otimes
(\Lambda^\bullet T y^{-1}q^{k-\frac ji})\otimes
(Sym^\bullet T^*q^{k-1+\frac ji})
$$
$$
\otimes
(Sym^\bullet T q^{k-\frac ji})
\bigr]
\bigr]
\bigr]^{\ZZ/i\ZZ}
\Bigr]
$$
$$
=p^iy^{-\dim X/2} \frac 1i \sum_{r=0}^{i-1}
\int_X \prod_{l=1}^{\dim X} x_l \prod_{k\geq 1}
\prod_{m=0}^{i-1} 
\frac
{(1-yq^{k-1+\frac mi}\xi^{mr}\ee^{-x_l})(1-y^{-1}q^{k-\frac mi}
\xi^{-mr}\ee^{x_l})}
{(1-q^{k-1+\frac mi}\xi^{mr}\ee^{-x_l})
(1-q^{k-\frac mi}\xi^{-mr}\ee^{x_l})}
$$
$$
=p^iy^{-\dim X/2}\frac 1i \sum_{r=0}^{i-1}
\int_X \prod_{l=1}^{\dim X} x_l
\prod_{j\geq 1} 
\frac
{(1-yq^{\frac {j-1}i}\xi^{(j-1)r}\ee^{-x_l})
(1-y^{-1}q^{\frac ji}\xi^{jr}\ee^{x_l})} 
{
(1-q^{\frac{j-1}i}\xi^{(j-1)r}\ee^{-x_l})
(1-q^{\frac ji}\xi^{jr}\ee^{x_l})
}
$$
$$
=p^i\frac 1i \sum_{r=0}^{i-1} Ell(X;y,q^{\frac 1i} \xi^r) = 
\sum_{m,l}c(mi,l)y^lq^m.
$$
Here we have denoted the primitive $i$-th root of unity by $\xi$.
Now Lemma \ref{DMVVlemma} 
finishes the proof. \end{proof}

\begin{remark}In \cite{wang} the authors conjectured an equivariant
version of \ref{DMVVtheorem}. Its proof follows 
using the same arguments as above. More precisely, we have the following.
Let $X$ and $G$ be as above and let $G \wr \Sigma_n$ be the wreath product 
(consisting of pairs $((g_1,...,g_n);\sigma), g_i \in G, \sigma \in \Sigma_n$ 
with multiplication:
$((g_1,...,g_n);\sigma_1) \cdot ((h_1,...,h_n);\sigma_2)=
((g_1 \cdot h_{\sigma_1^{-1}(1)},...,g_n \cdot h_{\sigma_1^{-1}(n)});
\sigma_1\sigma_2 )$). $G \wr \Sigma_n$ acts in an obvious way on $X^n$ and 
if $Ell_{orb}(X,G;y,q)=\sum c_G(m,l)y^lq^m$ then
\begin{equation}
\sum_{n\geq 0}p^n Ell_{orb}(X^n,G \wr \Sigma_n;y,q)=
\prod_{i=1}^\infty \prod_{l,m}
\frac 1
{(1-p^iy^lq^m)^{c_G(mi,l)}}.
\label{DMVVequivar}
\end{equation}
To obtain a proof of this formula, one should make the following 
changes in the above proof of Theorem \ref{DMVVtheorem}.
Using the description of the 
conjugacy classes in wreath products (cf. for example \cite{kerber})
$ \sum_{n\geq 0}p^n Ell_{orb}(X^n,G \wr \Sigma_n;y,q)$ can be rewritten as
the right hand side of the first row of \eqref{firststep}
with summation taken over collections $\{h \},a_1,...,a_n$ where 
$a_i$ as earlier are positive integers and $\{ h \}$ runs through all
 conjugacy classes in $G$. The same transformation as was used in 
\eqref{firststep} now yields the product over $i$ and $\{ h \}$ of terms
in which invariants are taken for the semidirect product 
of the centralizer of $h$ and ${\bf Z} /i{\bf Z}$ with the sheaf $V$
constructed for $X^h$. Finally, each term in this product is 
the graded dimension of a supersymmetric
algebra, which Lemma \ref{DMVVlemma} expresses in terms of 
$\chi_{i,\{ h \}}$. Calculation similar to the above calculation of 
$\chi_i$ identifies $\chi_{i,\{ h \}}$ 
with $\sum_{m,l} c_{\{ h \}}(mi,l)y^lq^m=
  y^{-dimX /2+F(h,X^h\subseteq X)} \frac 1 {|C(h)|}
\sum_{g\in C(h)} L(g, V_{h,X^h\subseteq X})
$ (the component of the orbifold elliptic genus corresponding to 
the conjugacy class $\{ h \}$). This yields \eqref{DMVVequivar}.
\end{remark}

\section{Comparison of different notions of elliptic genera}\label{comp.gen}
It is natural to ask how the orbifold elliptic genus of $X/G$ 
is related to its singular elliptic genus. To begin with, even in 
the case $|G|=1$, these two genera differ by a normalization factor.
In addition, when $\mu:X\to X/G$ has a ramification divisor 
$D=\sum_i(\nu_i-1)D_i$, one has to compare $Ell_{orb}(X,G;y,q)$ 
not to $\widehat{Ell}(X/G;y,q)$ but rather to 
$\widehat{Ell}(X/G,\Delta_{X/G};y,q)$ where 
$$\Delta_{X/G}:=\sum_j\left(\frac {\nu_j-1}{\nu_j}\right)\mu(D_j)$$
with the sum taken among representatives $D_j$ of the orbits of the action
of $G$ on the components of the ramification divisor.
\begin{conjecture}\label{mainconj}
{
Let $X$ be a smooth algebraic variety equipped with an effective action 
of a finite group $G$. Then
$$
Ell_{orb}(X,G;y,q)=
\left(
\frac 
{2\pi\ii\,\theta(-z,\tau)}
{\theta'(0,\tau)}
\right)^{\dim X}
\widehat{Ell}(X/G,\Delta_{X/G};y,q)
$$
where $\Delta_{X/G}$ is defined above.
}
\end{conjecture}

We will now present some evidence to support this conjecture.
\begin{proposition}
{ 
Conjecture \ref{mainconj} holds in the limit $\tau\to \ii\infty$.
}
\end{proposition}

\begin{proof} 
At $q=0$ the function $Ell_{orb}$ specializes to
$E_{orb}(y,1)$ of \cite{Batyrev}. Then the result of \cite{denefloeser} 
allows one to rewrite it in terms of $E_{st}(y,1)$, and Proposition 
\ref{normcross} finishes the proof.
\end{proof}

\begin{proposition}\label{mainconj.toric}
{
Conjecture \ref{mainconj} holds in the case when $X$ is a smooth toric 
variety and $G$ is a subgroup of the big torus of $X$.
}
\end{proposition}

\begin{proof}
Let $\Sigma$ be the defining cone of $X$ in the lattice $N$, see 
for example \cite{Danilov}. Let $n_i$ be the generators of one-dimensional 
cones of $\Sigma$. The group $G$ can be identified with $N'/N$ where
$N'$ is a suplattice of $N$ of finite coindex. Then the variety 
$X/G$ is given by the same cone $\Sigma$ in the new lattice $N'$.
The map $\mu:X\to X/G$ has ramification if and only if for some one-dimensional
rays of $\Sigma$ points $n_i$ are no longer minimal in the new lattice.

Torus-invariant divisors on a toric variety correspond to piecewise 
linear functions on the fan. It is easy to see that
the definition of $\Delta_{X/G}$ assures that the 
piece-wise linear linear function
that takes values $(-1)$ on all $n_i$ gives the divisor
$K_{X/G}+\Delta_{X/G}$. We denote this piece-wise linear function by 
$\deg$. One can show that 
$$Ell_{orb}(X,G;y,q)=
\left(
\frac 
{2\pi\ii\theta(-z,\tau)}
{\theta'(0,\tau)}
\right)^{\dim X}
f_{N',\deg z}(q)$$
where $f_{N',\deg z}(q)$ is the function defined in \cite{BorGun}.
More explicitly,
$$f_{N',\deg z}(q)=\sum_{m\in (N')^*} 
\Bigl(\sum_{C\in\Sigma}(-1)^{\codim C}
{\rm a.c.}{\sum_{n\in C\cap N'} q^{m\cdot n} \ee^{2\pi\ii z\deg(n)}}\Bigr),
$$
where $a.c.$ means analytic continuation.
The proof of this fact is based on the explicit calculation of the 
Euler characteristics of the bundles $V_{X^h\subseteq X}$ by means
of \v{C}ech cohomology. The calculation is very similar to that
of Theorem 3.4 of \cite{BorGun} and is left to the reader.
We remark that the sum over $h$ in Definition \ref{orbifold.def}
facilitates the change from $N$ to $N'$, while taking
$C(h)$-invariants is responsible for the switch from $N^*$
to its sublattice $(N')^*$.

Now let $Y\to X/G$ be a toric desingularization of $X/G$ given
by the subdivision $\Sigma_1$ of $\Sigma$. We denote the codimension
one strata of $Y$ by $E_k$, and the generators of the corresponding 
one-dimensional cones of $\Sigma_1$ by $r_k$. We also denote 
the first Chern classes of the corresponding divisors by $e_k$ and get
$$
\widehat{Ell}(X/G,\Delta_{X/G};y,q)=
\int_Y
\Bigl(\prod_l
\frac{(\frac{y_l}{2\pi\ii})\theta(\frac{y_l}{2\pi\ii}-z)\theta'(0)}
{\theta(-z)\theta(\frac{y_l}{2\pi\ii})}
\Bigr)\times
\Bigl(\prod_k
\frac{\theta(\frac{e_k}{2\pi\ii}-(\alpha_k+1)z)\theta(-z)}
{\theta(\frac{e_k}{2\pi\ii}-z)\theta(-(\alpha_k+1)z)}
\Bigr)
$$
where $c(TY)=\prod_l(1+y_l)$ and $\alpha_k=\deg(r_k)-1$.
We use $c(TY)=\prod_k (1+e_k)$ to rewrite 
$\widehat{Ell}(X/G,\Delta_{X/G};y,q)$ as
$$
\int_{Y}
\Bigl(\prod_k
\frac{(\frac{e_k}{2\pi\ii})\theta(\frac{e_k}{2\pi\ii}-\deg(r_k)z)\theta'(0)}
{\theta(-\deg(r_k)z)\theta(\frac{e_k}{2\pi\ii})}
\Bigr)
$$
which equals $f_{N',\deg z}(q)$ by Theorem 3.4 of \cite{BorGun}. We have 
used here the fact that $f$ does not change when the fan is subdivided.
\end{proof}

\begin{remark}
{
The above proposition was the main motivation behind our definition of 
the singular elliptic genus. 
}
\end{remark}

\begin{proposition}\label{mainconj.curve}
{
Conjecture \ref{mainconj} holds for $\dim X=1$.
}
\end{proposition}

\begin{proof}
Expanding $\theta$ functions as (linear) polynomials in
cohomology classes, one obtains that singular genus is equal to
$(2g-2)\theta'(-z)/(2\pi\ii\theta(-z))$ plus sum of contributions 
of singular points that depend on the ramification numbers only.
Here $g$ is the genus of $X/G$.
For the orbifold genus, one needs to notice that $h={\bf id}$ term gives
$(2g-2)\theta'(-z)/(2\pi\ii\theta(-z))$ plus contributions of points,
because it is the Euler characteristics of the bundle on the quotient
that is the usual elliptic genus bundle twisted at the ramification
points. 
Since the equality holds in the toric
case of the $d$-fold covering of $\PP^1$ by $\PP^1$, 
which has two points 
of ramification $(d-1)$, the extra terms for two genera coincide,
which finishes the proof.
\end{proof}

One would also want to compare singular elliptic genus to the 
elliptic genus defined for toric varieties and 
Calabi-Yau hypersurfaces in toric varieties in \cite{borlibg}. 
It turns out that these definitions agree, up to a normalization.
We will explain the Calabi-Yau case in more detail, and will leave
the toric case to the reader. We  need to recall the combinatorial
description of Calabi-Yau hypersurfaces in toric varieties and 
the previous definition of their elliptic genera.

Let $M_1$ and $N_1$ be dual free abelian groups of rank $d+1$. Denote
by $M$ and $N$ the dual free abelian groups $M=M_1\oplus {\bf
Z}$ and $N=N_1\oplus{\bf Z}$. Element $(0,1)\in M$ is denoted by $\rm deg$
and element $(0,1)\in N$ is denoted by $\deg^*$. There are dual reflexive
polytopes $\Delta\in M_1$ and $\Delta^*\in N_1$ which give rise to dual
cones $K\subset M$ and $K^*\subset N$. Namely, $K$ is a cone over
$(\Delta,1)$ with vertex at $(0,0)_M$, and similarly for $K^*$. 
There is a complete fan $\Sigma_1$ on $N_1$ whose one-dimensional cones 
are generated by some lattice points in $\Delta^*$ (in particular, by all 
vertices). This fan induces the decomposition of the cone $K^*$ into 
subcones, each of which includes $\deg^*$. Let us denote this decomposition 
by $\Sigma$. A generic Calabi-Yau hypersurface $X_f$ of the family 
given by the above combinatorial data is determined by a choice of 
coefficients $f_m$ for each $m\in (\Delta,1)$.

Elliptic genus of $X_f$ was defined in \cite{borlibg} as the graded Euler
characteristic of a certain sheaf of vertex algebras on $X_f$. We will
not need to recall the definition of this sheaf in view of the following 
combinatorial formula for the elliptic genus: 
\begin{proposition}
{ 
The elliptic genus $Ell(X_f;y,q)$ of the Calabi-Yau hypersurface $X_f$
as defined in \cite{borlibg} is given by
$$
Ell(X_f;y,q)=y^{-d/2}\sum_{m\in M} 
a.c.\left(
\sum_{n\in K^*}
y^{n\cdot {\deg} -m\cdot {\deg}^*}
q^{m\cdot n+m\cdot{\deg}^*}
G(y,q)^{d+2}
\right),
$$
where $a.c.$ stands for analytic continuation and 
$$G(y,q)=\prod_{k\geq 1}
\frac
{(1-yq^{k-1})(1-y^{-1}q^k)}{(1-q^k)^2}.
$$
}
\end{proposition}

\begin{proof} Combine Proposition 4.2 
and Definition 5.1 of \cite{borlibg}. \end{proof}

\begin{theorem}
Elliptic genus of the Calabi-Yau hypersurface $X_f$ of dimension $d$
defined above and its singular elliptic genus are related by the 
formula
$$
Ell(X_f;y,q)=
\left(
\frac 
{2\pi\ii\theta(-z,\tau)}
{\theta'(0,\tau)}
\right)^{d}
\widehat{Ell}(X_f;y,q).
$$
\end{theorem}

\begin{proof}
First of all, observe that 
$$y^{-1/2}G(y,q)=
\frac{2\pi\ii\theta(-z,\tau)}
{\theta'(0,\tau)},$$
due to the product formulas for $\theta(z,\tau)$ and $\theta'(0,\tau)$,
see \cite{Chandra}. Therefore, we only need to show that
$$
\widehat{Ell}(X_f;y,q)=
\sum_{m\in M} 
a.c.\left(
\sum_{n\in K^*}
y^{n\cdot {\deg} -m\cdot {\deg}^*}
q^{m\cdot n+m\cdot{\deg}^*}
G(y,q)^{2}
\right).
$$
Denote by $\deg_1$ the piece-wise linear function on $N_1$ whose 
value on the generators of the one-dimensional cones of $\Sigma_1$ is
$1$. Notice that $K^*$ consists of all points $(n_1,l)\in N$ such 
that $l\geq \deg_1(l)$. In addition, one can replace $\sum_{n\in K}\ldots$
by $\sum_{C\in\Sigma}(-1)^{\codim \Sigma}\ldots$ to get
$$
\sum_{m\in M} 
a.c.\left(
\sum_{n\in K^*}
y^{n\cdot {\deg} -m\cdot {\deg}^*}
q^{m\cdot n+m\cdot{\deg}^*}
G(y,q)^{2}
\right)
$$
$$
=
\sum_{k\in\ZZ}\sum_{m_1\in M}\sum_{C_1\in \Sigma_1}
(-1)^{\codim C_1}a.c.\sum_{n_1\in C_1}\sum_{l\geq \deg_1(n_1)}
y^{l-k}
q^{m_1\cdot n_1+lk+k}
G(y,q)^{2}
$$
$$
=
\sum_{k\in\ZZ}\sum_{m_1\in M}\sum_{C_1\in \Sigma_1}
(-1)^{\codim C_1}a.c.\sum_{n_1\in C_1}
\sum_{l\geq \deg_1(n_1)}
y^{\deg_1(n_1)-k}
q^{m_1\cdot n_1 + \deg_1(n_1)k+k}(1-yq^{k})^{-1}
G(y,q)^{2}
$$
$$
=G(y,q)^{2}
\sum_{k\in \ZZ}
\frac
{y^{-k}q^{k}}
{(1-yq^k)}
f_{N_1,\deg_1 z}(yq^k,q).
$$
Let $\Sigma_1'$ be a refinement of the fan $\Sigma_1$ in $N_1$ 
such that the corresponding toric variety $\PP_{\Sigma_1'}$ is smooth. 
Coefficients $f_m$ define a hypersurface $X_f'$ in $\PP_{\Sigma_1'}$
which is a resolution of singularities $X_f$. We denote the codimension
one strata of $\PP_{\Sigma_1'}$ by $D_j$, their first Chern classes by 
$d_j$ and the corresponding generators of one-dimensional cones of
$\Sigma_1'$ by $n_j$. By Theorem 3.4 of \cite{BorGun}, we get
$$G(y,q)^{2}
\sum_{k\in \ZZ}
\frac
{y^{-k}q^{k}}
{(1-yq^k)}
f_{N_1,\deg_1 z}(yq^k,q)
$$
$$
=G(y,q)^{2}
\sum_{k\in \ZZ}
\frac
{y^{-k}q^{k}}
{(1-yq^k)}
\int_{\PP_{\Sigma_1'}}
\prod_j
\frac{(\frac{d_j}{2\pi\ii})\theta(\frac{d_j}{2\pi\ii}-\deg_1(n_j)(z+k\tau))
\theta'(0)}
{\theta(-\deg_1(n_j)(z+k\tau))\theta(\frac{n_j}{2\pi\ii})}
$$
$$
=\int_{\PP_{\Sigma_1'}}
\prod_j
\frac{(\frac{d_j}{2\pi\ii})\theta(\frac{d_j}{2\pi\ii}-\deg_1(n_j)z)
\theta'(0)}
{\theta(-\deg_1(n_j)z)\theta(\frac{n_j}{2\pi\ii})}
\Bigl(
\sum_{k\in \ZZ} G(y,q)^2\frac{y^{-k}q^k}{(1-yq^k)}
\ee^{k\sum_j d_j\deg_1(n_j)}
\Bigr)
$$
$$
=\int_{\PP_{\Sigma_1'}}
\prod_j
\frac{(\frac{d_j}{2\pi\ii})\theta(\frac{d_j}{2\pi\ii}-\deg_1(n_j)z)
\theta'(0)}
{\theta(-\deg_1(n_j)z)\theta(\frac{n_j}{2\pi\ii})}
\Bigl(
\sum_{k\in \ZZ} G(y,q)^2\frac{y^{-k}q^k}{(1-yq^k)}
\ee^{k\sum_j d_j\deg_1(n_j)}
\Bigr)
$$
We denote $D=\sum_j \deg_1(n_1)D_j$ and $d=c_1(D)$.
Because of Proposition 3.2 of \cite{borlibg},
we get 
$$
\sum_{k\in \ZZ} G(y,q)^2\frac{y^{-k}q^k}{(1-yq^k)}
\ee^{k\sum_j d_j\deg_1(n_j)}=
\frac{G(e^dq,q)G(y,q)}{G(y^{-1}e^dq,q)}=
\frac 
{
2\pi\ii \theta(\frac d{2\pi\ii})\theta(-z,\tau)
}
{
\theta(\frac d{2\pi\ii}-z) \theta'(0)
}
$$
which gives
$$\sum_{m\in M} 
a.c.\left(
\sum_{n\in K^*}
y^{n\cdot {\deg} -m\cdot {\deg}^*}
q^{m\cdot n+m\cdot{\deg}^*}
G(y,q)^{2}
\right)
$$
$$=
\int_{\PP_{\Sigma_1'}}
\prod_j
\frac{(\frac{d_j}{2\pi\ii})\theta(\frac{d_j}{2\pi\ii}-\deg_1(n_j)z)
\theta'(0)}
{\theta(-\deg_1(n_j)z)\theta(\frac{n_j}{2\pi\ii})}
\Bigl(
\frac 
{
2\pi\ii \theta(\frac d{2\pi\ii})\theta(-z,\tau)
}
{
\theta(\frac d{2\pi\ii}-z) \theta'(0)
}
\Bigr)
$$
Observe now that $D=\pi^*(-K_{\PP_{\Sigma_1}})$
where $\pi:\PP_{\Sigma_1'}\to\PP_{\Sigma_1}$ is the resolution
induced by the subdivision of the fan. In addition, $X_f$ 
is a zero set of a section of $D$. Hence, the adjunction
formula gives 
$$c(T_{X'_f})=i^*c(T\PP_{\Sigma'_1})/(1+i^*d)$$
where $i:X'_f\to\PP_{\Sigma_1'}$ is the embedding. The exceptional divisors of 
$X_f'\to X_f$ are $D_j\cap X_f'$ (unless $\dim\pi(D_j)=0$),
and their discrepancies are
equal to $\deg(n_j)-1$. Then it is easy to see that the above expression
is precisely the singular elliptic genus of $X_f$.
\end{proof}

\begin{remark}
{ 
The case of toric varieties is a straightforward application of 
Theorem 3.4 of \cite{BorGun} and is left to the reader.}
\end{remark}

\begin{remark}
{
The above calculations indicate that for any smooth variety $\PP$ 
of dimension $d+1$
one can define a weak Jacobi form of weight $d$ and index $0$ which 
coincides with the singular elliptic genus of the Calabi-Yau
hypersurface in $\PP$ if $\PP$ has smooth anticanonical divisors.
Otherwise, the formula gives elliptic genus of ``virtual'' 
Calabi-Yau hypersurface in $\PP$. One can also interpret this Jacobi
form as an elliptic genus of $(d+1,1)$-dimensional Calabi-Yau
supermanifold $\Pi KX$ (canonical line bundle over $X$, considered 
as an odd bundle).}
\end{remark}

\section{Cobordism invariance of orbifold elliptic genus}
\label{cobordism.section}
We shall view $Ell_{sing}(X/G)$ and $Ell_{orb}(X,G)$ as invariants 
of $G$-action on $X$ and will work in the category of 
stably almost complex manifolds.

\begin{lemma} 
{
Singular elliptic genus is an  invariant of 
complex $G$-cobordism.
}
\end{lemma}

\begin{proof} We shall consider cobordisms of pairs $(X,D)$
(\cite{wall})
where $X$ is a stably almost complex manifold 
(i.e. $C^{\infty}$ manifold such that a direct sum of a trivial bundle
$\epsilon$  
with the differentiable tangent bundle $T_X$ admits a complex structure)
and $D=\cup D_i$ is a finite 
union of codimension one stably almost complex submanifolds 
(i.e. $T_{D_i} \oplus \epsilon$ is a complex subbundle in 
$\epsilon \oplus T_X$) satisfying the following normal crossing 
condition: at each 
point of $D_{i_1} \cap ... \cap D_{i_k}$ the union of (stabilized by 
adding trivial bundles)  
tangent spaces $T_{D_{i_j}} \oplus \epsilon$ 
is given in the (stabilized) tangent space to $X$ by $l_1 \cdot \cdot \cdot l_k=0$ where 
$l_i$ are linearly independent complex linear forms. 
A pair $(X,D)$ is cobordant to zero if
there exist a $C^{\infty}$-manifold $Y$ with a complex structure on 
the stable tangent bundle and a system of submanifolds $\cup E_i$ 
such that $\partial Y=X$ and $\cup \partial E_i=\cup D_i$. 
As usual, the disjoint union and product provide the ring structure 
on cobordism classes. 

Notice that the numbers $c_{i_1} \cup ... \cup c_{i_k} \cup 
[D_{j_1}] \cup ... \cup [D_{j_k}] ([X])$ where $[D_i]$  are the 
classes in $H^2(X,{\bf Z})$ dual to submanifolds $D_i$,  
$[X] \in H_{2\dim_{\bf C} \ X}$  is the fundamental class of $X$ 
and $\sum_j i_j +k=\dim_{\bf C} \ X$ are invariants of cobordism 
of such pairs (indeed, if $X =\partial Y$ and $j: X \rightarrow Y$ then 
this number is $j^*(c_{i_1} \cup ... \cup c_{i_k} \cup 
[E_{j_1}] \cup ... \cup [E_{j_k}]) ([X])=
c_{i_1} \cup ... \cup c_{i_k} \cup [D_{j_1}] \cup ... \cup [D_{j_k}] (j_*[X])
=0$
since $X$ is homologous to zero in $Y$).  The lemma therefore will follow 
if we shall show that for an almost complex null-cobordant $G$-manifold $X$ 
the quotient $X/G$ admits a resolution of singularities $(\tilde {X/G},D)$,
where $D=\cup D_i$ is the exceptional locus, such that this pair is cobordant 
to zero.

If $X=\partial Y$, where $Y$ is a $G$-manifold, we can construct resolution of 
$Y/G$ as follows. Let $H$ be a subgroup of $G$ and $Y_H=\{ y \in Y \vert Stab \ y =H \}$.
Then $Y_H$ are smooth submanifolds of $Y$ (possibly with boundary) 
providing a stratification of $Y$.    
Let $C(H)$ be the union of subgroups of $G$ conjugate to $H$. Then 
$Y_{C(H)}=\cup_{H' \in C(H)} Y_{H'}$ is still a submanifold of $Y$ and the 
group $G$ acts on $Y_{C(H)}$ so that $Y_{C(H)} \rightarrow Y_{C(H)}/G$ 
is an unramified cover (of degree $[G:H]$). In particular, $Y_{C(H)}/G$ 
is a smooth manifold and these manifolds for all $H \subset G$ provide
a stratification of $Y/G$ such that $Y/G$ is equisingular along each stratum 
$Y_{C(H)}/G$. A small regular neighborhood of each stratum in $Y/G$
is isomorphic to a bundle $\xi_H$ over $Y_{C(H)}/G$ with the fiber 
isomorphic to 
$V/H$ where $V$ is a fiber of the normal bundle to $Y_H$ in  $Y$ over a point
of $Y_H$ (this presentation is independent of a point in $Y_H$ and 
representations
at points of $Y_H$ and $Y_{H'}$ are isomorphic for conjugate $H$ and $H'$).
  
For each quotient singularity $V/H$ let us fix the universal desingularization
constructed by Bierstone-Milman (cf. Theorem 13.2 of \cite{bierstone}).
Its universality 
assures that it is equivariant
with respect to the centralizer of $H$ in $GL(V)$. 
Hence one can use the transition functions of $\xi_H$ to construct 
the fibration $\tilde \xi_H$ with the same base as $\xi_H$ and having as 
its fiber the universal resolution 
of $V/H$. Moreover, due to universality of canonical resolution
(cf. Theorem 13.2 in \cite{bierstone}) this property assures  
that $\tilde \xi_H$ corresponding to different 
classes of conjugate subgroups $H$ can be glued together yielding an 
almost complex manifold which boundary is the pair $(\tilde {X/G},D)$ 
where $D$ is the exceptional set of the universal resolution of $X/G$. 
This proves the lemma.
\end{proof}

\begin{lemma} 
{
Orbifold elliptic genus is an invariant of $G$-cobordism.
}
\end{lemma}

\begin{proof} Let $X$ be a null-cobordant $G$-manifold. Then for 
each $g \in G$ the  pair $X^g, \nu (X^g, X)$ where $\nu (X^g, X)$ is the 
normal bundle of the fixed point set $X^g$ in $X$ is cobordant to zero as
well. Since the contribution of the term in $Ell_{orb}$ corresponding to a conjugacy 
class $[g]$ is a combination of the products of Chern classes of $X^g$ and 
$\nu (X^g, X)$ evaluated on the fundamental class of $X^g$ this contribution is
zero. This yields the lemma.
\end{proof}

\begin{corollary} 
{
Conjecture \ref{mainconj} is true for $G={\bf Z}/2{\bf Z}$.
}
\end{corollary}

\begin{proof} This follows from the result of Kosniowski 
(\cite{kosn})
describing generators of ${\bf Z}/p{\bf Z}$-cobordisms.
If $p=2$, then additive generators of cobordism group
in any dimension are toric varieties with group
being a subgroup of the big torus. 
Hence Proposition \ref{mainconj.toric} yields the claim.
\end{proof}

\end{document}